\documentclass[12pt]{article}
\usepackage{latexsym}
\usepackage{amsfonts}
\usepackage{amsmath}
\usepackage{amssymb}

\newtheorem{thm}{Theorem}
\newtheorem{lem}[thm]{Lemma}
\newtheorem{cnj}[thm]{Conjecture}
\newtheorem{prop}[thm]{Proposition}
\newtheorem{cor}[thm]{Corollary}
\newtheorem{obs}[thm]{Observation}
\newtheorem{rem}[thm]{Remark}

\newcommand{\R}{\mathbb{R}}
\newcommand{\Rn}{\R^n}

\def\eps{\varepsilon}

\def\qed{\hfill $\vcenter{\hrule height .3mm
\hbox {\vrule width .3mm height 2.1mm \kern 2mm
\vrule width .3mm height 2.1mm} \hrule height .3mm}$ \bigskip}
\def\wrt{with respect to }

\begin{document}
\title{Duality of Metric Entropy
\thanks{ This research was
partially supported by grants from the US-Israel BSF (all authors)
and the NSF [U.S.A.] (the third-named author).
    }}
\date{}
\author{ S. Artstein,
V. Milman, S. J. Szarek
}
\maketitle
%
%
%
\abstract{For two convex bodies $K$ and $T$ in $\Rn$, the covering number
of $K$ by $T$,
denoted $N(K,T)$,
is defined as the minimal number of translates of $T$ needed to
cover $K$.
Let us denote by $K^{\circ}$ the polar body of
$K$ and by $D$ the euclidean unit ball in $\Rn$.
We prove that the two functions
of $t$, $N(K,tD)$ and $N(D, tK^{\circ})$, are equivalent in the
appropriate sense, uniformly over symmetric convex bodies
$K \subset \Rn$ and over $n \in \mathbb{N}$.
In particular, this  verifies the duality conjecture for entropy numbers
of linear operators,
posed by Pietsch in 1972, in the central case when either the domain or
the range of the operator
is a Hilbert space.
}

\begin{section}{Introduction}

For two convex bodies $K$ and $T$ in $\Rn$, the covering number
of $K$ by $T$,
denoted $N(K,T)$,
is defined as the minimal number of translates of $T$ needed to
cover $K$
\[ N(K,T) = \min\{N : \exists \, x_1 \ldots x_N \in \Rn , ~K\subset
\bigcup_{i\le N} \ x_i +T\}. \]
We denote by $D$ the euclidean unit ball in $\Rn$.
In this paper we
prove the following duality result for covering numbers.

\begin{thm}[Main theorem] \label{mainthm}
There exist two universal constants $\alpha$ and $\beta$ such that for
any dimension $n$ and any convex body $K \subset \Rn$, symmetric \wrt the
origin, one has
\begin{equation} \label{main}
       {N(D, {\alpha}^{-1}K^{\circ})}^{1\over  \beta} \le {N(K, D)} \le
{N(D, {\alpha}K^{\circ})}^{\beta}
\end{equation}
where $K^{\circ} := \{u\in \Rn : \sup_{x\in K} \langle x, u \rangle \le
1\}$ is the polar body of $K$.
\end{thm}
\noindent The best constant $\beta$ that our approach
yields is $\beta = 2+\eps$ for
any $\eps >0$,  with $\alpha  = \alpha (\eps)$.

Our theorem establishes a strong connection between
the geometry of a set and its polar or, equivalently, between a
normed space and its dual. Notice that since the theorem is true for any
$K$, we can actually infer that for any $t>0$
\begin{equation} \label{mainlog}
{\beta}^{-1}\log{N(D, {\alpha}^{-1}{t}K^{\circ})}
\le \log{N(K,t D)}
\le {\beta}\log{N(D, {\alpha}t K^{\circ})}.
\end{equation}
(For definiteness,
above and in what follows all logarithms are to the base~2.)
The quantity $\log{N(K,t \,T)}$ has a clear information-theoretic
interpretation: it is the complexity of $K$, measured in bits,
at the level of resolution $t$ \wrt the metric associated with $T$
(e. g., euclidean if $T=D$). Accordingly, (\ref{mainlog}) means that
the complexity of $K$ in the euclidean sense is controlled by
that of the euclidean ball \wrt $\|\cdot\|_{K^\circ}$
(the gauge of ${K^\circ}$, i.e., the norm whose unit ball is
${K^\circ}$), and vice versa, at {\em every } level of resolution. While
it is clear that these  complexities should be related, the universality
of the link that  we establish is somewhat surprising.

In addition to the immediate information-theoretic ramifications,
covering numbers appear in many other areas of mathematics. For
example, both  quantities $N(K, tD)$ and $N(D, tK^{\circ})$ enter the
theory of  Gaussian processes (see, e.g., \cite{dudley} and \cite{KL},
or the  survey \cite{ledoux}) and our results transform
some conditional statements into theorems (see, e.g., \cite{LL}).

Theorem \ref{mainthm} resolves an
old problem, going back to Pietsch (\cite{P}, p. 38) and
referred to as the ``duality conjecture for entropy
numbers," in a special yet most important case. The
problem
can be stated in terms of covering numbers in the
following way (below and in what follows we shall abbreviate
``symmetric \wrt the origin" to just ``symmetric").

\begin {cnj} [Duality Conjecture] \label{dualcnj}
Do there exist two numerical constants $a, b\ge 1$
such that for any dimension $n$,
and for any two symmetric convex bodies $K$ and $T$ in $\Rn$
one has
\begin{equation} \label{dc}
     \log{N(T^{\circ}, a K^{\circ})} \le b \log{N(K,T)} ,
\end{equation}
where $A^{\circ}$ denotes the polar body of $A$\,?
\end{cnj}
Theorem 1 verifies this conjecture in the case where one of
the  two bodies is a euclidean ball or, more generally, by affine
invariance of the problem, when one of the two bodies is an
ellipsoid. In the special case where {\em both} bodies
are ellipsoids it is well known and
easy to check that there is equality
in (\ref{dc}), with $a= b = 1$.

\smallskip
This conjecture originated in operator theory, and so
we restate it below in the language of
entropy numbers of operators. For two Banach spaces $X$ and $Y$,
with unit balls $B(X)$ and $B(Y)$ respectively,
and for a linear operator $u:X\to Y$, the $k^{th}$ entropy number of $u$
is defined by
$$
e_k(u) := \inf \{\eps : N(uB(X), \eps B(Y))\le 2^{k-1} \} .
$$
(In fact, above and in what follows $k$ does not need
to be an integer.)
So, for example, $e_1(u) = \|u\|_{op}$ (the operator norm),
and one can easily see that
$e_k(u) \to 0$ as $k \to \infty$ if and only if $u$ is a compact operator.
    Therefore the two sequences
$\left(e_k(u)\right)$ and $\left(e_k(u^*)\right)$
always begin with the same number
$\|u\|_{op} = \|u^*\|_{op}$, and
$e_k(u) \to 0$ if and only if $e_k(u^*)\to 0$.
Since the sequence $\left(e_k(u)\right)$ can be thought of as
quantifying the compactness of the
operator $u$,
it seems natural to ask to what extent do $\left(e_k(u)\right)$ and
$\left(e_k(u^*)\right)$  behave similarly.
This is the context in which the duality conjecture was originally
formulated,  and it read as follows.

\vskip 6pt
\noindent {\bf (Duality Conjecture in the language of entropy numbers) }
{\it Do there exist numerical constants $a,b\ge 1$,
such that for any
two Banach spaces $X$ and $Y$, any linear operator $u:X\to Y$ and any natural
number $k$, one has }
$$
e_{bk}(u^*)\le ae_k(u) \ ?
$$

\smallskip \noindent
The two formulations are equivalent since
considering the entropy numbers of
the dual operator $u^*: Y^* \to X^*$
means covering $(B(Y))^{\circ}$ with  (translates of) $\eps
(B(X))^{\circ}$, and since one can restrict oneself to
bodies which are convex hulls of a finite number of points
and thus lie in a finite dimensional space.
This is formulated explicitly in
Observation \ref{fewext} of  Section 2 below.
In other words, Theorem \ref{main}
verifies the duality conjecture (when expressed in terms of entropy
numbers) for the case in which one of the two spaces,
either the domain or the range of the operator $u$, is a Hilbert space.

{Some special cases of the problem have been studied before, and some
particular results were established, see, e.g.,
\cite{Art}, \cite{AMSz}, \cite{BPST}, \cite{GKS}, \cite{KM}, \cite{MSz},
\cite{MSz1}, \cite{PT}, \cite{Pi}, \cite{PiBook}, \cite{S},
\cite{T}.
We mention two of the above which
have special relevance to our approach: firstly
\cite{KM}, which shows the duality for entropy numbers when the
rank of the operator is (at most) comparable to the logarithm of the
covering  number, and secondly \cite{T},
which demonstrates a form of duality
involving some measures of the size
of {\em entire } sequences $(e_k(u)), (e_k(u^*))$.}

The proof of the theorem consists of three parts. The first part
is based on a fact
has already been formulated and proven in the required form in our paper
\cite{AMSz}, in which
we establish duality up to some factor $\gamma $ depending on the body
$K$. Next, this step is iterated, each time applied to a different
body (for example, a multiple of $K$ intersected with a euclidean
ball of some radius), and
we bound the covering number by a product of covering numbers of polars.
In the third and last step we shrink this product to a product of only
two or three factors, establishing duality with absolute constants.
Since this is a two sided inequality,
almost every statement is divided into two parts. However, there is
generally no
interplay between the two arguments, and the proofs of
the two sides of the inequality can be
read independently.

We wish to point out that different iterations
could be used. One of them is outlined in a short
note \cite{AMS1}. We use here the one that yields the
best constant $\beta$ in the exponent and may potentially
lead to a result that is optimal in that regard.

The paper is organized as follows. In section 2 we show
how duality is established up to constants depending on
the diameter of the set.
In sections 3 and 4 we first present an iterating scheme
which yields a bound for the covering number
in the form of a long product, and
then a telescoping argument that shrinks the
product  to a mere product of two terms. This
will complete the proof of the  main theorem.
Section 5 consists of various additions to
the proof.
First, we show how to  improve the constant $\beta = 6$,
given by the method described in sections 2-5,
to  $\beta = 2+\eps$  (for any $\eps>0$, and
with $\alpha = \alpha (\eps)$). Next,
we state a related conjecture and several results
associated with that conjecture.

\noindent{\it Remark on notation:} Unless
otherwise stated,
above and in what follows all constant appearing are universal
(notably independent of the dimension and of the particular convex body
or the operator that is being considered). If a constant
$c$ depends on some parameter $\theta$, we
will indicate that by writing
$c(\theta)$.

\end{section}

\begin{section}{A first step toward duality}

For a symmetric convex body $K\subset \Rn$, denote
by $k$ the logarithm of its covering number (so that
$N(K,D) = 2^k$), and define the parameter
\[\gamma(K) := \max \{ 1, M^*(K \cap D) \sqrt {n\over k}\},\]
where, as usual, $ M^*(A)$ denotes half the mean width of the set $A$,
that is, $M^*(A) = \int_{S^{n-1}}\sup_{y\in A} \langle u,y
\rangle \ d\mu(u)$, with
$\mu$ the normalized Haar measure on the sphere $S^{n-1}$.

The first step of the proof of the main theorem
is a duality result involving the parameter $\gamma$
instead of a universal constant $\alpha$.
The following lemma is a combination
of two statements, the first of which appeared in \cite{MSz}
and the second in
\cite{AMSz}.

\begin{lem}[First step]\label{wlemma}
There exists a universal constant $c_{2}>0$ and, for every $\eps >0$,
a constant $C_2(\eps)> 0$ such that,
for any dimension $n$ and for any symmetric convex body
$K \subset \Rn$, denoting $\gamma = \gamma (K)$, we have
\begin{equation}\label{wlem1}
N(K, D)\le N(D, {c_2\over \gamma}K^{\circ})^{3}
\end{equation}
and
\begin{equation}\label{wlem3}
N(D, {C_2(\eps)\gamma}K^{\circ})\le N(K, D)^{1+\eps}.
\end{equation}
\end{lem}

For a general body $K\subset \Rn$ the parameter
$\gamma$ can be as large as
$\sqrt{n\over k}$.
In Observation \ref{fewext} below we explain why
we can restrict
our considerations to a certain special class of convex bodies, namely
the convex hulls of not too many points. In the rest of the
section we will show that in this class there
are good  bounds for $\gamma$.

\begin{obs}\label{fewext}
For any convex body $K$, any set $S\subset K \subset S+D$
of cardinality
$N(K, D)$ and for any $\rho>0$ we have,
denoting $T = {\rm conv}(S)$,
\[N(D, (2\rho +2)K^{\circ} ) \le N(D, \rho T^{\circ}).\]
Similarly, if $N(K, D) > 1$, then we can find $S\subset K$ of cardinality
$N(K,D)$ such that ${\rm diam} (S) = {\rm diam} (K)$
and that, denoting $T = {\rm conv}(S)$, we have
\[ N(K, D) \le N(T, {1\over 2}D).\]

\end{obs}

\noindent {\it Remark:~} The argument does not require that the
original body lies
in a  finite dimensional
space (whereas a convex hull of finitely many points obviously does).
In particular, this shows the equivalence of the operator theoretic
formulation of the duality conjecture and the finite dimensional analogue
with universal constants.

\vskip 6pt

\noindent{\it Proof~}
Obviously
$T\subset K\subset T + D$. Denoting
$N(D, \rho T^{\circ}) = N$, we can pick
a $\rho$$T^{\circ}$-net $\{y_i\}_{i=1}^{N}$ for $D$,
i.e., $D\subset \cup_{i=1}^{N} y_i+ \rho T^{\circ}$.
We want to pass to a net {\em inside} $D$, for this notice that
$y_i + \rho T^{\circ}$ intersects $D$, say at a point $z_i$,
and that $\{z_i\}_{i=1}^{N}$ is a  $2\rho T^{\circ}$-net
for $D$.
We
claim that  $\{z_i\}_{i=1}^{N}$ is
a $(2\rho+2) K^{\circ}$-net of $D$.
Indeed,
for every $y\in D$ there exists a $z_i$ in the net such that
$y-z_i \in 2\rho T^{\circ}$, i.e.,
$\sup_{x\in T}(y-z_i, x) \le 2\rho$. Hence
(using that $y-z_i$ is in $2D$) we have
$\sup_{x\in T+D}(y-z_i, x) \le 2\rho+2$, which means precisely that
$\|y-z_i\|_{(T+D)^\circ}\le 2\rho +2$.
In particular, since $K\subset T+D$,
we see that $\|y-y_i\|_{K^\circ}\le \rho +2$, as required.
We conclude that
$N(D, (2\rho+2)K^{\circ})\le N$, and this verifies the first
part of the observation.

For the second part, we denote this time
$N =  N(K,D) $ and pick a 1-separated set
$\{ x_i \}_{i=1}^{N}$ in $K$ which realizes the diameter. We do this
simply by
choosing  two points, the distance between which is the diameter of $K$,
and completing them to a 1-separated set of cardinality $N$.
Again, this is possible since a maximal separated
set has at least as many elements as the minimal covering. Denote
$T = {\rm conv}\{x_i\}$. Since $\{x_i\}$
were $1$-separated,  $N(T,{1\over 2}D)\ge N$.
This completes the demonstration of the observation. $\hfill \square$

\vskip 6pt

The following proposition is an estimate for $\gamma(K)$
which is valid
whenever $K$
is  the  convex hull of $\le 2^k$ points in $RD$ and, in
addition,  has a covering number $\le 2^k$. It
was established in \cite{MSz}.
The general conjecture, which still remains open,
is that for this class of bodies the parameter $\gamma$
is bounded by a universal constant, regardless of the diameter
of the body. If this were true, Lemma \ref{wlemma} and
Observation \ref{fewext} would
imply the duality of entropy numbers (with $1+\eps$
in the exponent!).
We discuss the conjecture in Section 5; for
a more elaborate discussion and related results we refer the reader to
\cite{MSz}.

\begin{prop}[An $O(\log^3{R})$ estimate for $\gamma$]\label{logR}
There exists a universal constant $C_0$ such that if a set
$S\subset RD\subset \Rn$ (for some $R>1$) consists of $2^k$ points, and if
$N(K , D) \le 2^k$ for $K= {\rm conv}S$, then
\[M^*(K\cap D) \le C_0 (\log{R})^3 \sqrt{k\over n}.\]

\end{prop}

\vskip 6pt
Lemma \ref{wlemma}, Observation \ref{fewext} and
the above  Proposition
\ref{logR} can be combined as follows.
Denote $ \psi (x)= 2 C_2 (C_0\log^3{x}+1) + 2$, where $C_2=C_2(1)$
comes from Lemma \ref{wlemma}.

\begin{cor}[Duality up to $\psi(R)$]\label{dupsi}
If $K\subset RD\subset \Rn$ then
\begin{equation}\label{psi1}
N( D ,\psi(R) K^{\circ}) \le N(K,D)^2
\end{equation}
and
\begin{equation}\label{psi2}
N(K,D) \le N( D ,(1/\psi(R)) K^{\circ})^3.
\end{equation}
\end{cor}

\end{section}

\begin{section}{An iterating scheme}

In this section we present an iterating procedure that gives a  bound for
the covering number. The  first lemma is based on a
simple  geometric iteration  procedure
(and admits a variant which is valid in the
non-euclidean case; see Remark \ref{noneuc}).

\begin{lem}[Iterating procedure]\label{iter}
For any symmetric convex body $K\subset~\Rn$ and any
sequence $R_0 < R_1 < \cdots < R_s$,
\begin{equation}\label{iter1}
N(D, R_0 K^{\circ})\le N(D, R_{s}K^{\circ})\prod_{j=0}^{s-1}N(D, {R_j\over 2}
(K\cap {R_{j+1}}D)^{\circ}),
\end{equation}
and
\begin{equation}\label{iter2}
N(K, R_0D)\le N(K, R_sD) \prod_{j=0}^{s-1}N(2K\cap R_{j+1}D, R_jD).
\end{equation}
\end{lem}

\noindent{\it Proof~}
For (\ref{iter1}) consider the following inequality
\[N(D, R_0K^{\circ})\le N(D, {R_0\over 2}
{\rm conv}(K^{\circ}\cup {1\over R_1}D))
N({R_0\over 2} {\rm conv}(K^{\circ}\cup {1\over R_1}D),R_0K^{\circ}),\]
which follows from the
sub-multiplicativity of covering numbers: for every $A, B$ and
$C$ it is true that $N(A, B)\le N(A,C)N(C,B)$.
Rewriting the first term on the right hand side,
changing the convex hull in the second term to the Minkowski sum of sets
(which is bigger and thus harder to cover) and using  the rule
$N(A+C, B+C) \le N(A,B)$ leads to
\[N(D, R_0K^{\circ})\le  N(D, {R_0\over 2}
(K \cap {R_1}D)^{\circ}) N(D, R_1K^{\circ}). \]
Repeating the above argument another $(s-1)$ times yields (\ref{iter1}).

To show (\ref{iter2}) we first notice that
\[N(K, R_0D) \le N(K, R_1D)N(R_1D\cap 2K, R_0D), \]
where we use
the fact that
$ N(K, R_1D \cap 2K) = N(K, R_1D)$,
since the centers of a covering
by euclidean balls may always be assumed to lie inside $K$,
and also use
sub-multiplicativity
of covering numbers. Iterating  this inequality
gives (\ref{iter2}). The proof of Lemma
\ref{iter} is thus complete. $\hfill \square$
\vskip 12pt
Now is the time to choose the sequence $(R_j)$.
In fact, we will choose two different sequences, each corresponding
to a different inequality in the main theorem.
There is much freedom in this choice, and we do not
suggest that our choice is optimal.

For the first sequence, let $R_0$ be a large constant to be
specified later. Define $R_{j+1}$ by the
formula
\[{\sqrt{R_j}\over 2} = \psi\left({R_{j+1}\over \sqrt{R_j}} \right).\]
Remembering that
$\psi(x) = 2 C_2(C_0(\log x)^3+1) +2 $,  the above means that
\[
R_{j+1} = \sqrt{R_j} \ \exp\left(  (({\sqrt{R_j}-4-4C_2})/4C_2C_0)^{1/3}
\right).\]
In particular,  if $R_0$ is large enough
then this sequence increases to
$\infty$.  (This is needed since we
will later use the fact that  $N(D, R_{j}K^{\circ}) =
1$ for $j$ large enough.) Corollary \ref{dupsi} together with
Lemma \ref{iter} imply now the following

\begin{cor}\label{iterr1}
With the above choice of the sequence $(R_j)$ we have, for every
symmetric convex body $K$,
\begin{equation}\label{iter11}
N(D, R_0 K^{\circ})\le N(D, R_{s}K^{\circ})\prod_{j=0}^{s-1}
N(K\cap {R_{j+1}}D, \sqrt{R_j}D)^2.
\end{equation}
\end{cor}

\noindent{\it Proof~}
To deduce Corollary \ref{iterr1} from Lemma \ref{iter}
we only need  to  explain the inequality
\[N(D, {R_j\over 2}(K\cap {R_{j+1}}D)^{\circ})\le
N(K\cap {R_{j+1}}D, \sqrt{R_j}D)^2.\]
To this end, rewrite
\begin{eqnarray*}
N\left(D, {R_j\over 2}(K\cap {R_{j+1}}D)^{\circ}\right) & = &
N\left(D, {\sqrt{R_j}\over 2}\left({K\over \sqrt{R_j}}\cap
{R_{j+1}\over \sqrt{R_j}}D\right)^{\circ}\right)\\
& = & N\left(D, \psi\left({R_{j+1}\over \sqrt{R_j}}\right)\left({K\over
\sqrt{R_j}}\cap  {R_{j+1}\over \sqrt{R_j}}D\right)^{\circ}\right)\\
& \le & N\left({K\over \sqrt{R_j}}\cap
{R_{j+1}\over \sqrt{R_j}}D, D\right)^2 \\
& = & N(K\cap {R_{j+1}}D, \sqrt{R_j}D)^2,
\end{eqnarray*}
where for the inequality we used
(\ref{psi1}) of Corollary \ref{dupsi}. $\hfill \square$
\vskip 6pt
For the proof of the other side of the inequality in the main theorem
we have a different condition on the sequence to make
this type of argument work. Again, let
$R_0'$ be a big constant to be specified later.
Define  $R_{j+1}'$ by
\[ \psi\left({R_{j+1}'\over R_j'}\right) = {\sqrt{R_j'}\over 2},\]
which can be rewritten as $R_{j+1}' = R_j' \,
   \exp\left(  \left(\left({\sqrt{R_j'}-4-4C_2}\right)/4C_2C_0\right)^{1/3}
\right)$. Again, it is clear that this sequence is
increasing   to $\infty$.

\begin{cor}\label{iterr2}
With the above choice of a sequence $R_j'$ we have,
for every convex symmetric body $K$,
\begin{equation}\label{iter22}
N(K, R_0' D)\le N(K, R_{s}'D)\prod_{j=0}^{s-1}
N(D, \sqrt{R_j'}(K\cap {R_{j+1}'\over 2}D)^{\circ})^3.
\end{equation}
\end{cor}

\noindent{\it Proof~}
Again, we will use Lemma \ref{iter} together with Corollary \ref{dupsi}.
Here we should explain the inequality
\[N(2K\cap R_{j+1}'D, R_j'D)\le
N(D, \sqrt{R_j'}(K\cap {R_{j+1}'\over 2}D)^{\circ})^3.\]
This is even simpler since
\begin{eqnarray*}
N(2K\cap R_{j+1}'D, R_j'D)
&=& N\left({2\over R_j'}K\cap {R_{j+1}'\over R_j'}D, D\right) \\
&\le& N\left(D, {R_j'\over 2\psi\left({R_{j+1}'\over R_j'}\right)}
(K\cap {R_{j+1}'\over 2} D)^{\circ}\right)^3\\
& = & N(D, \sqrt{R_j'}
(K\cap {R_{j+1}'\over 2} D)^{\circ})^3,
\end{eqnarray*}
where for the inequality we used
(\ref{psi2}) of Corollary \ref{dupsi}. $\hfill \square$

\end{section}

\begin{section}{Telescoping the long product}  \label{telescoping}

In this last step we collapse the long products of covering numbers
appearing in (\ref{iter11}) and (\ref{iter22}) to
products consisting of just two terms. The largest $R_s$
(respectively $R_s'$) will be chosen to
exceed the diameter of
the set, and so the terms $N(K, R_s'D)$ and $N(D, R_sK^{\circ})$
will  both equal 1.
We need the following two super-multiplicativity inequalities
for covering numbers which are valid for any symmetric convex body $K$.

\begin{lem}\label{telesc}
Let $A>a>3B>3b$. Then
\begin{eqnarray}
N(K\cap AD, aD) N(K\cap BD, bD) &\le& N(K\cap AD, {b\over 4}D)
\label{triang1}\\
\label{triang2}
N(D, a(K\cap AD)^{\circ}) N(D, b(K\cap BD)^{\circ}) &\le &
N(D, {b\over 4}(K\cap AD)^{\circ}) .
\end{eqnarray}
\end{lem}
\noindent{\it Proof~}
Since $K$ enters the inequlities only 
via its intersections with balls 
of radii $\le A$, we may as well assume that 
$K = K\cap AD$ to begin with. 
For the first inequality, denote $N_1 = N(K, aD)$ and
$N_2 = N(K\cap BD, bD)$. Pick an $a$-separated set $x_1, \ldots x_{N_1}$
in $K$ and a $b$-separated set $y_1, \ldots y_{N_2}$
in $K\cap BD$ (both separations with respect to the euclidean norm).
Define a new set by $z_{i,j} = x_i/2 +y_j/2$.
All these points are in $K$, and there are $N_1N_2$ of them.
We shall show that, in addition,  $z_{i,j}$'s are $({b/ 2})$-separated;
this will imply $N(K, {b\over 4}D) \ge N_1N_2$, as required.
To show the asserted separation, we consider two cases. First, if we look
at
$|z_{i,j}-z_{i,k}|$, this is simply $|y_j-y_k|/2$ and it exceeds
$b/2$. On the other hand, if   $k\neq i$, then
$|z_{i,j}-z_{k,l}| \ge |x_i - x_k|/2 - |y_j-y_l|/2$, and using the fact
that the $y_i$'s are in $BD$ we see that these quantities are greater than
$a/2 - B$, which in turn exceeds
${b\over 2}$.
This completes the proof of
inequality (\ref{triang1}).

For the second inequality in the Lemma, denote
$N_1 = N(D, aK^{\circ})$ and $N_2 = N(D, b(K\cap BD)^{\circ})$.
Pick sets $\{x_1, \ldots, x_{N_1}\}$ and $\{y_1, \ldots y_{N_2}\}$
in $D$ which are respectively
${a}K^{\circ}$-separated and
${b}(\alpha K^{\circ}+{1-\alpha \over B}D)$-separated,
where $\alpha = { a\over 2a - b} \in (0,1)$
(note that  $\alpha K^{\circ}+{1-\alpha \over B}D
\subset {\rm conv}(K^{\circ}\cup
{1 \over B}D) = (K \cap BD)^\circ$).
Define $z_{i,j} = {b\over 2a}x_i + (1-{b\over 2a})y_j$. All
these points are in $D$,
and there are  $N_1N_2$ of them. As above, it will be enough to show that
the $z_{i,j}$'s are $ {b\over 2}K^{\circ}$-separated,
i.e., whenever
$(i,j) \neq (k,l)$, then
\[z_{k,l} \not\in z_{i,j} + {b\over 2}K^{\circ}. \]
When looking at $j=l$, this is
the same as asking that
${b\over 2a}x_k \not\in {b\over 2a} x_i + {b\over 2}K^{\circ}$,
which follows from the separation of
$x_i$ and $x_k$.
When looking at  $j\neq l$ and noticing that $x_i, x_k \in D$,
we see that it suffices to show that
\[ (1-{b\over 2a})y_j \not\in (1-{b\over 2a})y_l + 2{b\over 2a}D +
    {b\over 2}K^{\circ}. \]
Under our hypotheses, the above follows from the
separation of $y_l$ and $y_j$. Indeed, ${1 \over 2}(1-{b\over 2a})^{-1}
=
\alpha$ just by the definition of $\alpha$. On the other hand, it is
readily verified that the assumption
$a>3B>2B+b$ implies ${1 \over a}(1-{b\over 2a})^{-1} < {{1 - \alpha} \over
B}$. The proof is thus
complete. $\hfill \square$

\begin{rem}\label{noneuc}~
(i)
With more careful argument, any factor less than $1/2$ instead
of ${1/4}$ can be obtained in the
lemma (with then stronger conditions on $a,B$).
Moreover, if we work  with separated sets instead of covering numbers,
then we may arrive at any factor less than 1: for a factor
$1-\eps$ (which corresponds to ${1\over 2} - {\eps \over 2} $
for covering numbers), we need the condition $a> 3({1-\eps\over \eps})B$.
This may be used to  improve slightly the constants in our
main theorem.

(ii) Notice that $D$ plays no special role
in Lemma \ref{telesc};
the same inequalities hold for two general symmetric
convex bodies $K$ and $T$ (i.e., with
$D$ replaced by an arbitrary $T$).
\end{rem}
\vskip 12pt

\noindent{\it Proof of the Main Theorem~}
We will successively apply Lemma \ref{telesc}
to the long products in (\ref{iter11}) and (\ref{iter22}).
However, an additional trick is required since for two neighboring
factors in the products the
condition of Lemma \ref{telesc} does not hold, and so they cannot be
``collapsed."
For example, for two such factors in (\ref{iter11})
one has $a = \sqrt{R_j}$
and $B= R_j$, and so one cannot hope for $a>3B$.
The trick is to split the product into two parts,
by grouping separately the factors corresponding to the odd and the
even $j$'s. The growth of $R_j$ is fast enough
so that the conditions of Lemma \ref{telesc} are satisfied
for each two consecutive odd factors, and
for each two consecutive even factors.
We provide details for the product from (\ref{iter11});
the analysis of  (\ref{iter22}) is fully analogous.

First choose $s$ to be the smallest even number so that $R_s > {\rm diam}(K)$.
Then the product in (\ref{iter11}) which bounds $N(D, R_0K^{\circ})$ can
be written as (we omit the power $2$ for the moment)
\[ \prod_{j=1}^{s/2} N(K\cap R_{2j}D, \sqrt{R_{2j-1}}D)
\prod_{j=1}^{s/2} N(K\cap R_{2j-1}D, \sqrt{R_{2j-2}}D).\]
For the first collapsing step in each of the two sub-products we need to
check that $\sqrt{R_{s-1}}> 3R_{s-2}$, and that $\sqrt{R_{s-2}}>
3R_{s-3}$.  When
$R_{j+1} = \sqrt{R_j}\,
   \exp\left(  (({\sqrt{R_j}-4-4C_2})/4C_2C_0)^{1/3}
\right)$,
these conditions clearly hold as
long as $R_{s-3}$ is larger than some numerical constant $C$.
Since $R_j > R_0$ for $j > 0$, it is enough to start with
$R_0$  which is big enough.
(To be able to later compare the
obtained expressions with $N(K,D)$, we also insist that $R_0 \ge 16$.)
After this first step, using  Lemma \ref{telesc},
the product becomes (bounded by)
\begin{eqnarray*}
&&N(K\cap R_sD, {\sqrt {R_{s-3}}\over 4}D)
\prod_{j=1}^{s/2-2} N(K\cap R_{2j}D, \sqrt{R_{2j-1}}D) \\
&&N(K\cap R_{s-1}D, {\sqrt {R_{s-4}}\over 4}D)
\prod_{j=1}^{s/2-2} N(K\cap R_{2j-1}D, \sqrt{R_{2j-2}}D).
\end{eqnarray*}
   From here onward all the steps are the same;
we just need to make sure at each stage
$j$ that
\begin{equation}\label{condi}
(\sqrt {R_j}/4)/(R_{j-1}/2)\ge 3.
\end{equation}
As before, this is indeed
satisfied if $R_j >C$, which is assured
since we insist that $R_0>C$.
We point out that the factors ${1/4} $
in $b/4$ do not accumulate, but enter into the
quantity $a$ of the next step.
Continuing this way with all the factors of these products,
we arrive at
\[ N(D, R_0K^{\circ}) \le N(K\cap R_sD, {\sqrt {R_0}\over 4}D)^2
N(K\cap R_{s-1}D, {\sqrt{R_1}\over 4}D)^2, \]
(we have inserted back the power $2$ in (\ref{iter11}))
which, having  insisted that $R_0 \ge 16$, implies
\[N(D, R_0K^{\circ}) \le N(K, D)^4.\]
Similarly, in the other direction we use Corollary \ref{iterr2}
and inequality (\ref{triang2}) to obtain
\[N(K, R_0D) \le N(D, K^{\circ})^6,\]
and the proof of the main theorem is complete. $\hfill \square$
\vskip 6pt
As mentioned earlier, a large part of this proof carries over to the
case of two {\em general} convex bodies. We summarize this in the
following conditional proposition.
(Our decision to include this statement in the form below was
influenced by discussions with Nicole Tomczak-Jaegermann.)

\begin{prop} Let $T$ be a convex symmetric body in a euclidean space such
that, for some constants $c, C>0$, the following holds:
       if $K$ is a convex symmetric body with $K\subset 4 T$, then
        \[ N(K, T) \le N(T^{\circ}, c K^{\circ})^{C}.\]
       Then, for some other constants $c',C' >0$ (depending only on $c, C$)
and \ {\bf any} convex symmetric body $K$
        \[ N(K, T) \le N(T^{\circ}, c'K^{\circ})^{C'}.\]
Dually, if $K$ is fixed and the hypothesis holds for all $T$'s verifying
$K\subset 4 T$, then the assertion holds for {\bf any}  $T$.
\end{prop}

\noindent{\it Proof~}
The argument consists of two parts. The first part is essentially a copy
of  the proof of Lemma \ref{iter} for the choice $R_{j} = 2^j$.
The only difference is that  an extra factor $2$ appears
in the analogue of (\ref{iter2}) since  at each step
\begin{eqnarray*}
N(K, 2^j T) &\le& N(K, 2^{j+2}T \cap 2K )N(2^{j+2}T\cap 2K, 2^jT)\\
&\le& N(K, 2^{j+1}T) N(2^{j+2}T\cap 2K, 2^jT)
\end{eqnarray*}
as we can no longer  assume that the centers of the covering are
inside $K$.
In this way we show the  inequality (similar to (\ref{iter2}))
\[ N(K, T) \le N(K,2^s T) \prod_{j=0}^{s-1}  N(2^{1-j}K\cap 4T, T)\]
and, dually, another inequality similar to (\ref{iter1}),
\[ N(K, T) \le N(K,2^s T)
       \prod_{j=0}^{s-1}
       N(K, {\rm conv}(2^{j-1}T \cup {1\over 4}K)).\]
For each factor in these products  the body that is being covered is
included  in $4$ times the covering body, and so we may use the assumption
(as before, we take $s$ to be the smallest integer such that
$N(K,2^s T)=1$) to pass to a product
of dual covering numbers. Thus for example we get that
\[ N(K, T) \le
       \prod_{j=0}^{s-1}
       N(T^{\circ}, c2^{j-1}(K\cap 2^{j+1}T)^{\circ}) ^C\]
(and a similar estimate if we use the second inequality instead).
We  then collapse the remaining product in the same way as in the
euclidean case, using Remark \ref{noneuc} (ii).
Note that we may now have to split the product into more
than 2, say
$l$, subproducts to make sure that all neighboring factors in each
subproduct satisfy  the condition of Lemma \ref{telesc}. However, the
resulting $l$  depends only on $c$. We thus arrive (in both cases) at
\[ N(K, T) \le
       N(T^{\circ}, {c\over 8} K^{\circ})^{C\log_2
       (48/c)}.\square \]


\end{section}

\begin{section}{{Improving the constant in the exponent}}

In this section we explain how to improve the constant $\beta = 6$
in (\ref{main}) that we have obtained in sections 2-4 to a constant
$\beta = 2(1+\eps)$.
The proof presented above
is somewhat non-symmetric.
As described, in one of the inequalities we get  $\beta = 6$,
and in the other $\beta = 4$. This $\beta = 4$ can be improved to
$2(1+\eps)$ if we continue working with a general  $\eps$ in
Lemma \ref{wlemma} and do not specify (as we did only to lighten the
notation) $\eps = 1$.
However, as stated above, we cannot in
a straightforward way obtain $\beta = 2(1+\eps)$ in the other
inequality. Below
we explain how to arrive at $\beta = 2(1+\eps)$ there
as well. We take this opportunity to elaborate upon the
conjecture called ``the Geometric Lemma," which we mentioned
in passing in Section~2.

\begin{cnj}[Geometric Lemma]\label{thegeol}
There exists an absolute constant $C_0$ such that for every
dimension $n$ and for every set $S= \{ x_i \}_{i=1}^{2^k}\subset
\Rn$ verifying $N(K, D) \le 2^k$,  where $K = {\rm conv}(S)$, we have
\[ M^*(K\cap D) \le C_0\sqrt{k\over n}.\]
\end{cnj}
\noindent Notice
that this is precisely a version of Proposition \ref{logR}
with no dependence on the radius of the set. We believe that its
importance may transcend its relevance to entropy numbers.

Considering now the dual situation, we can
formulate a dual version of the Geometric Lemma,
substituting the condition
$N(K, D) \le 2^k$ by the condition
$N(D, K^{\circ}) \le 2^k$.
However, having proved the duality of entropy numbers, it is easy to see
that the Geometric
Lemma and its dual version are formally equivalent.
Moreover, every estimate such as
Proposition \ref{logR} can be applied now to the dual situation.
For example the main theorem together with Proposition \ref{logR}
gives the following

\begin{prop}[A dual $O((\log R)^3)$ estimate for $M^*(K\cap
D)$]\label{dualogR} There exists a universal constant $C_0$ such that if a
set
$S\subset RD\subset \Rn$ consists of $2^k$ points, and if
$N(D, K^{\circ}) \le 2^k$
for $K= {\rm conv}S$, then
\[M^*(K\cap D) \le C_0 (\log{R})^3 \sqrt{k\over n}.\]
\end{prop}
\vskip 6pt
\noindent We can thus define a new parameter, $\gamma'(K)$, to be
\[\gamma'(K) := \max\{ 1, M^*(K\cap D) \sqrt{n\over \log N(D,
K^{\circ})}\},
\] and
repeating the argument of Lemma \ref{wlemma} we obtain
the following

\begin{lem}[Dual First step]\label{duwlemma}
        For every $\eps >0$
there is a constant $C_2'=C_2'(\eps)$ such that
for any dimension $n$ and for any symmetric convex body
$K \subset \Rn$, denoting $\gamma' = \gamma'(K)$ we have
\begin{equation}\label{duwlem3}
N(K, {C_2'\gamma'}D)\le N(D, K^{\circ})^{1+\eps}.
\end{equation}
\end{lem}

Employing the same line of argument as earlier, but using
Proposition \ref{dualogR} as an estimate on $\gamma'$, and
inequality  (\ref{duwlem3}) at every step, we are now able to obtain
$\beta = 2+\eps$, instead of $\beta = 6$, also in the
other inequality involved in the duality of metric entropy.

\vskip 12pt

To end this section, and the paper, we present another
proposition which is an application of both our results and our methods,
and gives an interesting link between Geometric Lemma type results and
the behavior of covering numbers under projections. Its proof
follows the same lines as that of Lemma
\ref{wlemma}, and other variants involving additional parameters
are possible.
Below we use the standard jargon of the asymptotic theory of normed
spaces, saying  that a property is satisfied for a ``random" projection
of given rank $k$ if it holds  for a set of (orthogonal)
rank $k$ projections whose
measure tends to $1$ as the relevant  parameters ($k, n$ below) tend to
infinity (where ``measure" =  ``the normalized Haar measure on the
corresponding Grassmann manifold").

\begin{prop} \label{lastprop}
There exist universal constants $c_{1},
C_{1}$, and for every $\lambda$ there exists a constant
$C_{2}=C_2(\lambda)$  depending only on $\lambda$ such that, for any
$K\subset \Rn$ with
$N(K,D) = 2^k$ and any integer $t_0$ with $k\le t_{0}\le n$ we have

(i) If  $M^*(K\cap D) \le \sqrt{t_{0}\over n}$,
then for every integer $t$ with
$t_{0} \le t \le n$, the random rank
$t$ projection $P_{t}$ satisfies
\[ N(P_{t}K, C_{1}\sqrt{t\over n}D)\le 2^k \,\,\,{\rm and}\,\,\,
N(P_{t}K, c_{1}\sqrt{t\over n}D)\ge 2^k.\]

(ii) In the other direction, if the random rank $t_{0}$
projection  $P_{t_{0}}$ verifies
\[ N(P_{t_{0}}K, \lambda \sqrt{t_{0}\over n}D)\le 2^k, \]
then necessarily
$M^*(K\cap D) \le C_{2}\sqrt{t_{0}\over n}$
and, for {\em any}  integer $t$ with $t_{0} \le  t \le n$, the random
rank $t$ projection $P_{t}$ satisfies
\[ N(P_{t}K, C_{2}\sqrt{t\over n}D)\le 2^k \,\,\,{\rm and}\,\,\,
N(P_{t}K, c_{1}\sqrt{t\over n}D)\ge 2^k.\]
\end{prop}

Thus we observe - as is typical in the asymptotic geometric analysis - a
unified form of  behavior for {\em all} dimensions and {\em all} convex
bodies.

We note that our Proposition \ref{logR} implies, in the case when
$K$ is a  convex hull of $2^k$ points, that
the critical $t_{0}$ in the above
proposition is bounded from above
by $C_{0}k(\log  k)^{6}$, for details see \cite{MSz} (to pass from
estimates on the  diameter to estimates on $\log N(K, D)$).
Notice that the validity of Conjecture \ref{thegeol} would imply that
for this class of bodies in fact $t_{0}\le C_{0}k$ for a universal $C_{0}$.
Also, our main theorem
implies that Proposition \ref{lastprop} remains true if we replace the
condition  on $N(K, D)$ with a similar one on $N(D, K^{\circ})$ (with
additional  universal constants). Similarly, we may replace the estimates
on the behavior of  covering numbers under projections with their dual
analogues, describing  the behavior of covering numbers under
intersections with random subspaces.

\end{section}

{\footnotesize
{S. Artstein, V. D. Milman

School of Mathematical Sciences,
Tel Aviv University,

Tel Aviv 69978, Israel

E-mail: artst@post.tau.ac.il

E-mail: milman@post.tau.ac.il

\vskip 8pt

S. J. Szarek

Equipe d'Analyse Fonctionnelle, B.C. 186,
Universit\'{e} Paris VI,

4 Place Jussieu,  F-75252  Paris, France

{\em and }

Department of Mathematics,
Case Western Reserve University,

Cleveland, OH 44106-7058, U.S.A.

E-mail: szarek@ccr.jussieu.fr
}}

\end{document}